\documentclass[12pt,oneside,english]{amsart}
\usepackage[latin1]{inputenc}
\usepackage{amssymb}

\makeatletter
\newtheorem{theorem}{Theorem}

\newtheorem{example}{Example}
\newtheorem{definition}{Definition}

\newtheorem{question}{Question}
\newtheorem{corollary}{Corollary}

\usepackage{babel}

\makeatother
\begin{document}
\baselineskip=17pt

\title[Binomial predictors]{Binomial predictors}

\author{Vladimir Shevelev}
\address{Departments of Mathematics \\Ben-Gurion University of the
 Negev\\Beer-Sheva 84105, Israel. e-mail:shevelev@bgu.ac.il}

\subjclass{11B37}

\begin{abstract}
 For a prime $p$ and nonnegative integers $n,k,$ consider the set $A_{n,\enskip k}^{(p)}=\{x\in [0,1,...,n]:\enskip p^k||\binom {n} {x}\}.$ Let the expansion of $n+1$ in base $p$ be: $n+1=\alpha_{0} p^{\nu}+\alpha_{1}p^{\nu-1}+...+\alpha_{\nu},$ where $0\leq \alpha_{i}\leq p-1,\enskip i=0,...,\nu.$ Then the number $n$ is called a \slshape  binomial predictor in base $p$\upshape ,\enskip if $|A_{n,\enskip k}^{(p)}|=\alpha_{k}p^{\nu-k},\enskip k=0,1,...,\nu.$ We give a full description of the binomial predictors in base $p.$
 \end{abstract}

\maketitle

\section{Introduction}
Let $p$ be a prime. For nonnegative integers $n,k,$ consider the set
 \begin{equation}\label{1}
 A_{n,\enskip k}^{(p)}=\{x\in [0,1,...,n]:\enskip p^k||\binom {n} {x}\}.
 \end{equation}
 Quite recently (2008), W.B. Everett [1] solved the following important problem.
\begin{question} How, knowing $n,$ to find the finite sequence
$$|A_{n,\enskip 0}^{(p)}|, |A_{n,\enskip 1}^{(p)}|, ... $$
(of course, without the direct calculations)?
  \end{question}
  Due to generality, the Everett's formula is sufficiently complicated.
  In this paper we indicate an infinite set of $n's,$ for which  Question 1 has especially simple solution which immediately follows from the expansion of $n+1$ in base $p.$ Conversely, in the limits of this set, knowing the sequence $\{|A_{n,\enskip k}^{(p)}|\},$ we can "predict" the expansion of $n+1$ in base $p.$ In the connection with this, we introduce the following notion.
\begin{definition}\label{1}
   Let the expansion of $n+1$ in base $p$ be: $n+1=\alpha_{0} p^{\nu}+\alpha_{1}p^{\nu-1}+...+\alpha_{\nu},$ where $0\leq \alpha_{i}\leq p-1,\enskip i=0,...,\nu.$ Then the number $n$ is called a \slshape  binomial predictor in base $p$\upshape ,\enskip if $|A_{n,\enskip k}^{(p)}|=\alpha_{k}p^{\nu-k},\enskip k=0,1,...,\nu.$
   \end{definition}
   \begin{example}\label{1}
   It is easy to see that $n=0$ is a binomial predictor in every base $p.$
   \end{example}
    Indeed, $|A_{0,\enskip0}^{(p)}|=1, $ that is the binary expansion of 1 in every base.
    \begin{example}\label{2}
   Let $p=2,\enskip n=11.$
   \end{example} Then $n+1=8+4.$
   The row of the binomial coefficients $\{\binom {n} {x}\enskip, x=0,1,...,11\}$ is:
   $$1, 11, 55, 165, 330, 462, 462, 330, 165, 55, 11, 1. $$
   Here $|A_{11,0}^{(2)}|=8, \enskip |A_{11,1}^{(2)}|=4, \enskip |A_{11,k}^{(2)}|=0,\enskip k\geq2.$ Thus, by the definition, 11 is a binomial predictor in base 2.\newline
   \begin{example}\label{3}
   Let $p=3,\enskip n=23.$ Then $n+1=2\cdot3^2+2\cdot3.$
    \end{example}
    Here $|A_{23,0}^{(3)}|=18, \enskip |A_{23,1}^{(3)}|=6, \enskip |A_{23,k}^{(3)}|=0,\enskip k\geq2.$ Thus 23 is a binomial predictor in base 3.\newline\newline
   Our aim is to give a full description of binomial predictors in base $p.$
    \begin{definition}\label{2}
A nonnegative integer $n$ is called a Zumkeller's number in base $p$ (in the case of $p=2$ see sequence $ A089633 \enskip in \enskip [12]$), if either it is 0 or its expansion in base $p$ has all digits $p-1,$ except, maybe, one; if the exceptional digit is the first ( it can occur only in case of $p\geq3$), then it could take an arbitrary value from $1,...,p-2;$ otherwise, it is only $p-2.$
   \end{definition}
Our result is the following.
\begin{theorem}\label{1}
$n\geq0$ is a binomial predictor in base $p$ if and only if it is a Zumkeller's number in the same base.
\end{theorem}
\section{Some classical results on binomial coefficients}
 The binomial coefficients play a very important role in numerous questions of number theory. For example, it is very known proof of the beautiful Chebyshev's theorem, using the binomial coefficients (see, for example, a Finsler's version of the proof in [13]). A connection between some questions of divisibility of the binomial coefficients and
 the old conjecture of the infinity of tween primes is appeared in the author's article [10].\newline
 The first important contributions into theory of binomial coefficients belong to Legendre (1830), Kummer (1852) and
 Lucas (1878). Let $p$ a prime and $a_p(n)$ be such exponent that
 \begin{equation}\label{2}
 p^{a_p(n)}||\enskip n.
  \end{equation}
  Let, furthermore,
 \begin{equation}\label{3}
n=n_0p^m+n_1p^{m-1}+...+n_m,\enskip0\leq t_i\leq p-1.
  \end{equation}
  be the expansion of $n$ in base $p.$ Denote\newpage
   \begin{equation}\label{4}
 s_p(n)=n_0+n_1+...+n_m,
  \end{equation}
  A.-M. Legendre [7, p.12] empirically noticed that (in our notations)
  \begin{equation}\label{5}
  a_p(n!)=(n-s_p(n))/(p-1).
  \end{equation}
  A proof see, e.g., in [9].
  From (5) we immediately obtain:
  $$ a_p(\binom {n} {x})=a_p(n!)-a_p(x!)-a_p((n-x)!)=$$
  $$((n-s_p(n))-(x-s_p(x))-(n-x-s_p(n-x))/(p-1)=$$
  \begin{equation}\label{6}
  (s_p(x)+s_p(n-x)-s_p(n))/(p-1).
  \end{equation}
   \begin{example}\label{4}
   \end{example}
  Since $s_2(2n)=s_2(n),$ then for the central binomial coefficients we find
  \begin{equation}\label{7}
  a_2(\binom {2n} {n})=s_2(n);\enskip a_2(\binom {2n+1} {n})=s_2(n+1)-1.
   \end{equation}
   Notice that in [11] was posed a question which remains open up to now.
   \begin{question}
  Does the diophantine equation $s_p(n)=s_q(n),$ where $p\neq q$  are fixed primes, have infinitely many solutions?
  \end{question}
  It easy to see that, by (6), the following equalities are equivalent:
  $$s_p((p-1)n)= s_q((q-1)m),$$
  $$ (p-1)a_p(\binom {pn} {p})=(q-1)a_q(\binom {qm} {q}).$$
  Furthermore, note that (6) implies the following simple corollary.
    \begin{corollary}\label{1}
 For every lattice pair $(x,y)\geq(0,0),$ we have the triangle inequality:
 $$ s_p(x+y)\leq s_p(x)+s_p(y).$$
 The equality attains if and only if $\binom {x+y} {x}$ is not multiple of $p.$
  \end{corollary}
   Now we can treat of Question 1 in a different foreshortening. Consider the equation
 \begin{equation}\label{8}
 s_p(x)+s_p(n-x)-s_p(n)=k(p-1), \enskip x\in[0,1,...,n].
 \end{equation}
 For $k\geq0, $ denote $\lambda_p^{(k)}(n)$ the number of solutions of (8). Thus we see that
 \begin{equation}\label{9}
|A_{n,\enskip k}^{(p)}|=\lambda_p^{(k)}(n),\enskip k=0,1,...
 \end{equation}
  In 1852, Kummer [6] made an important observation (a proof one can find, e.g., in [3]):
   \newpage \slshape$ a_p(\binom {n} {x})$ is the number of "carries" which appear in adding $x$ and $n-x$ in base $p.$\upshape \newline This statement plays a large role in the Everett's matrix method which was used for receiving
his general formula.\newline
 Another important result was obtained by Lucas [8]. He proved that if together with (3)
 \begin{equation}\label{10}
t=t_0p^m+t_1p^{m-1}+...+t_m,\enskip 0\leq t_i\leq p-1.
 \end{equation}
 then
 \begin{equation}\label{11}
\binom {n} {t}\equiv\prod_{i=0}^{m}\binom {n_i} {t_i}\pmod p.
 \end{equation}
 From (11) immediately follows the next corollary ([2]).
  \begin{corollary}\label{2}
 \begin{equation}\label{12}
\lambda_p^{(0)}(n)=(n_0+1)(n_1+1)...(n_m+1).
 \end{equation}
  \end{corollary}
  \bfseries Proof. \mdseries Indeed, (12) gives the number of all nonzero products in (11), when $0\leq t_i\leq n_i,\enskip i=0,...,m . $ Since $t_i, n_i\in [0,...,p-1],$ then none of the considered products is divisible of $p.\blacksquare$\newline
  Note that in the binary case in (12) sufficiently consider only factors with $n_i=1.$ Therefore, we have
   \begin{equation}\label{12}
\lambda_2^{(0)}(n)=2^{s_2(n)}.
 \end{equation}
 The latter is known result of J.Glaisher (1899; in [4] A.Granville gives a new elegant proof; generalizations in other direct see in [4]-[5]).

\section{Proof of necessity}
A proof of necessity we easily derive from \enskip\slshape only\upshape \enskip the first condition of Definition 1, i.e. from the equality $|A_{n,0}^{(p)}|=\alpha_0p^{\nu}.$  In connection with (3), note that $\nu\neq m$  only when $n_0=n_1=...=n_m=p-1.$ This a Zumkeller's number in base $p.$ In other cases we consider the equality $|A_{n,0}^{(p)}|=\alpha_0p^{m}.$ In its turn,$\alpha_0 \neq n_0\leq p-2$ only when $n_1=...=n_m=p-1.$ In this case we also have  a Zumkeller's number in base $p.$ Let now $n$ be a binomial predictor in base $p,$ when $n_0=p-1,$ but not all $n_i$ equal $p-1.$ Thus, we have
  \begin{equation}\label{13}
|A_{n,0}^{(p)}|=(p-1)p^{m}.
 \end{equation}
 From this equality and Corollary 2 we conclude that exactly $m$ from $m+1$ brackets in product (12) equal $p,$ while some one bracket equals $p-1.$ This means that exactly $m-1$ digits from $ n_1,..., n_m $ equal to $p-1,$ while some one digit equals $p-2.$ Thus in this case we again have  a Zumkeller's number in base $p.\blacksquare$\newpage

\section{Proof of sufficiency}
 Here we use the Kummer's theorem in the following equivalent form:\newline \slshape$ a_p(\binom {n} {x})$ is the number of "carries" which appear in subtracting $x$ from $n$ in base $p.$\upshape \newline
 Let $n$ be an Zumkeller's number. Evidently, in the trivial case when $n=p^{m+1}-1=\underbrace{(p-1)\vee (p-1)\vee...\vee (p-1)}_{m+1}$ we have a binomial predictor in base $p$ ( we use $\vee$ as operator of concatenation). Indeed, here $n+1=p^{m+1}$ and, by (9) and (12), $|A_{n,\enskip 0}^{(p)}|=p^{m+1}.$ On the another hand, for every $x\leq n,$ in the subtracting $x$ from $n$ in base $p$ we have no any "carries", i.e. $|A_{n,\enskip k}^{(p)}|=0,\enskip k\geq1.$  Analogously, we have a binomial predictor in case $n=n_0\vee \underbrace{(p-1)\vee (p-1)\vee...\vee (p-1)}_{m},\enskip n_0\leq p-2.$
 Indeed, here $n+1=(n_0+1)p^m $ and, by (9) and (12), $|A_{n,\enskip 0}^{(p)}|=(n_0+1)p^{m}.$ and again we have no any "carries" during the subtracting $x$ from $n$ in base $p.$ Let now
$n$ have a unique digit $p-1$ in its expansion in base $p.$ Consider such $n$ of the general form:
 \begin{equation}\label{15}
n=\underbrace{(p-1)\vee ...\vee (p-1)}_{t}\vee(p-2)\vee\underbrace{(p-1)\vee ...\vee (p-1)}_{m-t}.
\end{equation}
Then
$$n+1=\underbrace{(p-1)\vee ...\vee (p-1)}_{t+1}\vee \underbrace{0\vee ...\vee 0}_{m-t}=$$
\begin{equation}\label{16}
(p-1)(p^m+p^{m-1}+...+p^{m-t}).
\end{equation}
Let the $x$ in base $p$ has the form:
\begin{equation}\label{17}
x_0\vee ...\vee x_{t-1}\vee x_t\vee x_{t+1}\vee...\vee x_m.
\end{equation}
If $x_t\leq p-2,$ then in subtracting $x$ from $n$ in base $p$ we have not "carries". Evidently, that the number of such $x's$ is
\begin{equation}\label{18}
|A_{n,0}^{(p)}|=p^t(p-1)p^{m-t}=(p-1)p^{m}.
\end{equation}
If $x_t=p-1,$ such that also
\begin{equation}\label{19}
x_t=x_{t-1}=x_{t-2}=...=x_r=p-1, \enskip
\end{equation}
then, in view of $x\leq n,$ the length of this chain is not more than $t,$ i.e.
\begin{equation}\label{20}
1\leq r\leq t.
\end{equation}
In this case the number of "carries" is equal to the length of the chain, i.e. $t-r+1,$ and, by Kummer's theorem, we have
\begin{equation}\label{21}
x\in A_{n,t-r+1}^{(p)}.
\end{equation}\newpage
Putting
\begin{equation}\label{22}
t-r+1=k,
\end{equation}
we easily calculate $|A_{n,k}^{(p)}|:$
\begin{equation}\label{23}
|A_{n,k}^{(p)}|=p^{m+1-k}\cdot \frac {p-1} {p}=(p-1)p^{m-k}, \enskip k=1,...,t,
\end{equation}
where the factor $\frac {p-1} {p}$ corresponds to the digit $x_{r-1},$ i.e. the place after the end of chain (19).
Comparison of (18) and (23) with (16) shows that $n$ is a binomial predictor.
This completes the proof of Theorem 1. $\blacksquare$

\bfseries Acknowledgment.\mdseries \enskip The author is grateful to J.-P. Allouche for sending a scanned copy of
pp. 10-12 of Legendre book [7], and to D. Berend for important advises.

\end{document}